\documentclass[11pt]{article}
\usepackage{amsfonts}
\begin{document}
   \newtheorem{tw}{Theorem}[section]
   \newtheorem{lem}[tw]{Lemma}
   \newtheorem{wn}[tw]{Corollary}
   \newtheorem{uw}[tw]{Remark}
   \newtheorem{df}[tw]{Definition}
   \newtheorem{cl}{Claim}
   \newcommand{\bl}{\begin{lem}}
   \newcommand{\el}{\end{lem}}
   \newcommand{\bd}{\begin{df}}
   \newcommand{\ed}{\end{df}}
   \newcommand{\bu}{\begin{uw}\rm}
   \newcommand{\eu}{\end{uw}}
   \newcommand{\bc}{\begin{cl}\rm}
   \newcommand{\ec}{\end{cl}}
    \newcommand{\bnum}{\begin{enumerate}}
    \newcommand{\enum}{\end{enumerate}}
    \newcommand{\cg}{{\cal G}}
    \newcommand{\ct}{{\cal T}}
    \newcommand{\cgb}{{\cg_{\bullet}}}
    \newcommand{\sgb}{{\Sigma_{\bullet}}}
    \newcommand{\sfb}{{\Sigma^f_{\bullet}}}
    \newcommand{\x}{{\rm x}}
    \newcommand{\X}{{\rm X}}
    \newcommand{\y}{{\rm y}}
    \newcommand{\Y}{{\rm Y}}
    \newcommand{\z}{{\rm z}}
    \newcommand{\N}{{\mathbb N}}
    \newcommand{\Q}{{\mathbb Q}}
    \newcommand{\q}{{\bf Q}}
    \newcommand{\supp}{{\rm supp\:}}
    \newcommand{\spa}{{\rm span}}
    \newcommand{\spac}{\overline{\spa}}
    \newsymbol\varsubsetneq 2320
    \newsymbol\varsupsetneq 2321
    \newsymbol\square 1003

\begin{center}{\Large\bf Subsymmetric sequences and minimal spaces}
\end{center}

\begin{center}{\Large Anna Maria Pelczar}
\end{center}
\vspace{1cm}

\noindent {\small {\bf Summary}: We show that every Banach space
saturated with subsymmetric sequences contains a minimal
subspace.} \vspace{0.3cm}

\section{\large Introduction}

W.T.Gowers proved in \cite{dych} the celebrated dichotomy
concerning unconditional sequences and hereditarily indecomposable
spaces using Ramsey-type argumentation. In \cite{ram} he
generalized the reasoning and  showed, as an application, a
dichotomy concerning quasi-minimal spaces, ie. in which any two
infinitely dimensional subspaces contain further two infinitely
dimensional subspaces which are isomorphic. Putting these results
together he obtained the following "classification" theorem:

\begin{tw} {\rm \cite{ram}} Let $E$ be an infinitely dimensional
Banach space. Then $E$ has an infinitely dimensional subspace $G$
with one of the following properties, which are mutually exclusive
and all possible:

\vspace{0.2cm}

\noindent {\rm (1)} $G$ is a hereditarily indecomposable space,

\vspace{0.2cm}

\noindent {\rm (2)} $G$ has an unconditional basis and every isomorphism between
block subspaces of $G$ is a strictly singular perturbation of the
restriction of some invertible diagonal operator on $G$,

\vspace{0.2cm}

\noindent {\rm (3)} $G$ has an unconditional basis and is strictly
quasi-minimal (ie. does not contain a minimal subspace),

\vspace{0.2cm}

\noindent {\rm (4)} $G$ is a minimal space.\end{tw}

Natural question appears concerning the extensions of this
theorem. In this paper we prove that every Banach space saturated
with sub\-symme\-tric sequences contains a minimal subspace. It
follows that the class (3) can be restricted to strictly
quasi-minimal spaces not containing sub\-symme\-tric sequences and
it brings further division of the class (4) in terms of containing
sub\-symme\-tric sequences. An example of a minimal space not
containing sub\-symme\-tric sequences is the dual to Tsirelson's
space (\cite{lt}, \cite{cjt}), whereas Tsirelson's example is a
strictly quasi-minimal space (\cite{co}).

The method used here extends the technic applied in \cite{p},
which reflects the technic of Maurey's proof of Gowers' dichotomy
for unconditional sequences and HI spaces (\cite{dychm}). The same
method provides extensions also in the class (1) by examing
unconditional-like sequences introduced in \cite{t} (\cite{p}).

\vspace{0.5cm}

We introduce now the basic notation and definitions. Let $E$ be a
Banach space. Denote by $B_{E}$ the closed unit ball, by $S_{E}$ -
the unit sphere of $E$. Given a set $A\subset E$ by $\spa (A)$
(resp. $\spac(A)$) denote the vector subspace (resp. the closed
vector subspace) spanned by $A$. We will denote by $\Theta$ the
origin in the space $E$ in order to distinct it from the number
zero.

We say that two Banach spaces $E_1, E_2$ are $c-$isomorphic, for
$c\geq 1$, if there is an isomorphism $T: E_1\rightarrow E_2$
satisfying $\frac{1}{c}\| x\|\leq\| Tx\|\leq c\| x\|$ for $x\in
E_1$. Similarly we say that sequences $\{ x_n\}_n$, $\{ y_n\}_n$
of vectors of a Banach space are $c-$equiva\-lent, for $c\geq 1$,
if the mapping
$$T:\; \spa\{ x_n\}_n\ni \sum_{i=1}^n a_i x_i\rightarrow
\sum_{i=1}^n a_i y_i\in\spa\{ y_n\}_n$$ satisfies $\frac{1}{c}\|
x\|\leq\| Tx\|\leq c\| x\|$, $x\in E$.

Assume now that $E$ is a Banach space with a basis $\{ e_{n}\}_{n=1}^{\infty}$.

A support of a vector  $x=\sum_{n=1}^{\infty} x_{n} e_{n}$ is the
set $\supp x =\{ n\in\N :\; x_{n}\not =0\}$. We use notation $x<y$
for vectors $x,y\in E$, if every element of $\supp x$ is smaller
than every element of $\supp y$, $x<L$ for a vector $x\in E$ and a
subspace $L\subset E$, if every element of $\supp x$ is smaller
than every element of a support of any vector in $L$, and so forth
in this manner. A block sequence with respect to $\{ e_{n}\}$ is
any sequence of non-zero finitely supported vectors
$x_{1}<x_{2}<\dots$ , a block subspace - a closed subspace spanned
by a block sequence. We will use letters $x,y,z,\dots$ to denote
vectors of a Banach space, letters x,y,z,$\dots$ to denote finite
block sequences and capital letters X,Y,Z,$\dots$ for infinite
block sequences. Letters
$L,M,N,\dots$ will denote closed infinitely dimensional subspaces. For any finite
block sequence x by $|\x|$ denote the length of x, ie. the number
of elements of x. Given any two block sequences $\{ x_1,\dots , x_n\}<\{ y_1,y_2,\dots \}$ put
$$\{ x_1,\dots , x_n\}\cup\{ y_1,y_2,\dots \}=\{ x_1,\dots , x_n,y_1,y_2,\dots \}$$

For the convenience in the reasoning presented in the next section
we will treat $\{\Theta\}$ as a block sequence and adopt the
following convention: $|\Theta |=0$, $\Theta <x$ for any
$x\neq\Theta$, $\{\Theta\}\cup\{ y_1,y_2,\dots\}=\{
y_1,y_2,\dots\}$ for any block sequence $\{ y_1,y_2,\dots\}$.

We will work on a special class of block subspaces spanned by a
dense subset of $E$. By ${\bf Q} $ denote the vector space over
$\Q$, if $E$ is a real Banach space, or over $\Q+i\Q$, if $E$ is a
complex Banach space, spanned by the basis $\{ e_{n}\}_n$.
Obviously {\bf Q} is a countable dense set in $E$.

Denote by $\cg (E)$ the family of all infinitely dimensional and
closed subspaces of $E$. By $\cgb (E)$ denote the family of all
infinitely dimensional block subspaces spanned by block sequences
of vectors from the set ${\bf Q}$. Given a subspace $M\in\cgb (E)$
put $$\cgb (M)=\cgb (E)\cap\cg (M).$$

Given a subset $A\subset E$ let $\Sigma (A)$ (resp. $\Sigma^f
(A)$) be the set of all infinite (resp. finite) block sequences
contained in $A$. Put
$$\sgb (A)=\Sigma (A)\cap {\bf Q}^{\N},\;\;\;\sfb
(A)=\Sigma^f(A)\cap \bigcup_{n\in\N}{\bf Q}^n.$$
The family $\sgb (E)$ can be identified with the family $\cgb (E)$ in the obvious way.

While restricting our consideration to the family of block
sequences we will use a standard fact:

\bl\label{bbases} Let $E$ be a Banach space with a basis $\{
e_i\}_i$. Let $\{ x_n\}_n\subset E$ be a sequence satisfying
$\lim_{n\rightarrow\infty}e^*_i(x_n)=0$, $i\in\N$, where $\{
e^*_i\}_i$ is the sequence of  biorthogonal functionals of $\{
e_i\}_i$. Then for any $\varepsilon >0$ there is a block sequence
$\{ y_n\}_n$ which is $(1+\varepsilon)-$equivalent to some
subsequence of the sequence $\{ x_n\}_n$. \el

\section{\large The "stabilizing" Lemma}

In this section we present the key Lemma for our paper. It
reflects some combinatorial technics used in \cite{dychm},
\cite{dychz}.

Define a quasi-ordering relation on the family $\cg (E)$: for
subspaces $L,M\in\cg (E)$ write $L\leq M$ iff there exists a
finitely dimensional subspace $F$ of $E$ such that $L\subset M+F$.
This relation induces an equivalence relation: $L\doteq M$, for
$L,M\in\cg (E)$, if $L\leq M$ and $M\leq L$. In our consideration
we use a simple observation: for any subspaces $L,M\in\cg (E)$
satisfying $L\leq M$ we have $L\cap M\doteq L$, in particular
$L\cap M\in\cg (M)$.

We will prove now the Lemma, generalizing the argumentation given
in the proof of some properties of "zawada" (Lemma 1.21) in
\cite{dok}, which uses a standard now dia\-go\-na\-li\-za\-tion.

\bl\label{stab} {\rm \cite{p}} Let  $E$ be a Banach space with a basis. Let
$\tau$ be a mapping defined on the family $\cgb (E)$ with values in
the family $2^{\Sigma}$ of subsets of some countable set $\Sigma$.

If the mapping $\tau$ is monotone with regard to the relation
$\leq$ in $\cgb (E)$ and the inclusion $\subset$ in $2^{\Sigma}$,
i.e.
$$\forall N,M\in\cgb (E):\;\;\; N\leq M\Longrightarrow \tau (N)\subset \tau (M) $$
$${\rm or}\;\;\;\;\forall N,M\in\cgb (E):\;\;\; N\leq M\Longrightarrow \tau (N)\supset \tau (M) ,$$
 then there exists a subspace $M\in\cgb (E)$ which is stabilizing for $\tau$, i.e.
$$\forall L\in\cgb (M) :\;\;\; \tau (L)=\tau (M). $$\el

\noindent {\it Proof.} We can assume that the mapping $\tau$ is
increasing. If the mapping $\tau$ is decreasing, then put $\tau
'(M)=\Sigma\setminus\tau (M)$, a stabilizing subspace for $\tau '$
will be also stabilizing for $\tau$.

Suppose that for any subspace $N\in\cgb (E)$ there exists a further
subspace $L\in\cgb (N)$ such that $\tau (L)\varsubsetneq\tau (N)$. We
will construct a transfinite sequence $\{ L_{\xi}\} \subset \cgb
(E)$, indexed by the set of ordinal numbers $\{\xi :\; \xi
<\omega_{1}\}$, where $\omega_{1}$ is the first uncountable
ordinal, such that
$$\xi <\eta\;\Longrightarrow\; L_{\eta}\leq L_{\xi},\;\tau (L_{\eta})\varsubsetneq\tau (L_{\xi})$$

For $\xi =0$ put $L_{0}=E$. Take an ordinal number $\xi
<\omega_{1}$ and assume that we have defined subspaces $L_{\eta}$
for $\eta <\xi$. We consider two cases:

1. $\xi$ is of the form $\eta +1$. Then by our hypothesis there
exists a subspace $L_{\xi}\subset L_{\eta}$ such that $\tau
(L_{\xi})\varsubsetneq \tau (L_{\eta})$.

2. $\xi$ is a limit ordinal number. Since $\xi <\omega_{1}$, $\xi$
is a limit of some increasing sequence $\{\xi_{n}\}_n$ of ordinal
numbers (Theorem 5, 8.2 \cite{km}).

By the induction hypothesis we have
$$L_{\xi_{1}}\geq L_{\xi_{2}}\geq L_{\xi_{3}}\geq \dots\;\;\;\;\; {\rm and}\;\;\;\;\;
\tau (L_{\xi_{1}})\varsupsetneq \tau (L_{\xi_{2}})\varsupsetneq \tau (L_{\xi_{3}})\varsupsetneq \dots$$

By the monotonicity of the sequence $\{ L_{\xi_{n}}\}_n$ we have
$L_{\xi_{1}}\cap \dots\cap L_{\xi_{n}}\in\cgb (E)$ for $n\in\N$.
Choose a block sequence $\{ a_{n}\}_n$ such that $a_{n}\in
L_{\xi_{1}}\cap \dots\cap L_{\xi_{n}}\cap\q$ for $n\in\N$ and define
$L_{\xi}=\spac \{ a_{n}\}_{n\in\N}$. Then obviously $L_{\xi}\leq
L_{\xi_{n}}$ for $n\in\N$, hence $\tau
(L_{\xi})\subset\tau(L_{\xi_{n}})\varsubsetneq\tau
(L_{\xi_{n-1}})$, which ends the construction.

Hence we have constructed an uncountable family $\{ \tau
(L_{\xi})\}_{\xi <\omega_{1}}$ of strongly decreasing (with
respect to the inclusion) subsets of the set $\Sigma$, which
contradicts the countability of $\Sigma$.\hfill $\square$

\bu Let $E$ be a Banach space. Notice
that Lemma \ref{stab} holds also for the family of all block
subspaces or the family $\cg (E)$. \eu

\section{\large Sub\-symme\-tric sequences and minimal spaces}

\bd A Banach space $E$ is called $c-$minimal, for $c\geq 1$, if
any infinitely dimensional closed subspace of $E$ contains a
further subspace which is $c-$iso\-morphic to $E$.

A Banach space $E$ is called {\it minimal}, if any infinitely
dimensional closed subspace of $E$ contains a further subspace
which is isomorphic to $E$.\ed

\bd A basic sequence $\{ x_n\}_{n\in\N}\subset E$ is called
$c-$sub\-symme\-tric, for $c\geq 1$, if it is unconditional and is
$c-$equivalent to any of its subsequence.

A basic sequence $\{ x_n\}_{n\in\N}\subset E$ is called {\it
subsymmetric}, if it $c-$sub\-symme\-tric for some $c\geq 1$. \ed

While restricting to block sequences we will apply the following

\bl\label{bbasessub} Let $E$ be a Banach space with a basis. If
$E$ contains a $C-$sub\-symme\-tric sequence, for some constant
$C\geq 1$, then for any $\delta >0$ the space $E$ contains also a
$(C+\delta )-$sub\-symme\-tric block sequence.\el

\noindent {\it Proof of Lemma \ref{bbasessub}.} Let $\{ e_i\}_i$
be a basis for $E$. By $\{ e^*_i\}_i$ denote the biorthogonal
functionals for $\{ e_i\}_i$. Let $\{ x_n\}_n\subset E$ be a
$C-$sub\-symme\-tric sequence, for $C\geq 1$. We can assume,
picking a subsequence of $\{ x_n\}_n$ if needed by Cantor diagonal
method, that for some scalars $\{ a_i\}_i$ we have
$\lim_{n\rightarrow\infty}e^*_i(x_{n})=a_i$, $i\in\N$. Put
$z_n=x_{2n}-x_{2n-1}$ for $n\in\N$. Then $\{ z_n\}_n$ is a basic
unconditional sequence. Take any strictly increasing function
$\phi :\N\rightarrow\N$. Define a strictly increasing function
$\psi:\N\rightarrow\N$ as follows: for any $n\in\N$ put $\psi
(2n)=2\phi (n)$ and $\psi (2n-1)=2\phi (n)-1$. Notice that the
corresponding isomorphism $T:\spac\{ x_n\}_n\rightarrow\spac\{
x_{\psi (n)}\}_n$ given by subsymmetry of $\{ x_n\}_n$, satisfies
$$T(z_n)=T(x_{2n}-x_{2n-1})=x_{\psi (2n)}-x_{\psi (2n-1)}=
x_{2\phi (n)}-x_{2\phi (n)-1}=z_{\phi (n)},\;\;\; n\in\N$$
Hence the sequence $\{ z_n\}_n$ is also $C-$sub\-symme\-tric.

Fix $\delta >0$. Pick $\eta >0$ such that $(1+\eta )^2C<C+\delta$.
Obviously $\lim_{n\rightarrow\infty}e^*_i(z_n)=0$, $i\in\N$, hence
by Lemma \ref{bbases} there is a block sequence $\{ y_n\}_n$ which
is $(1+\eta )-$equivalent to some subsequence of $\{z_n\}_n$. Thus
by the choice of $\eta $ the sequence $\{ y_n\}_n$ is $(C+\delta
)-$sub\-symme\-tric. \hfill $\square$

\vspace{0.5cm}

We say that a Banach space is saturated with sequences of a given
type, if every its subspace contains a sequence of this type.

Now we present the main results:

\begin{tw}\label{sub} Let $E$ be a Banach space saturated with
$C-$sub\-symme\-tric sequences, for some $C\geq 1$. Then for any
$\varepsilon >0$, the space $E$ contains a $(C^2+\varepsilon
)-$minimal subspace. \end{tw}

\begin{wn} A Banach space saturated with sub\-symme\-tric sequences
contains a minimal space. \end{wn}

\noindent {\it Proof of Corollary.} Let $E$ be a Banach space
saturated with sub\-symme\-tric sequences. By the standard
diagonal argumentation there is a subspace $E_0\in\cg (E)$ which
is saturated with $C-$sub\-symme\-tric sequences for some $C\geq
1$. Indeed, if this was not the case, one could choose a
decreasing sequence of subspaces $\{ E_n\}_n\subset\cg (E)$ such
that  for $n\in\N$ the space $E_n$ contains no
$n-$sub\-symme\-tric sequence. Let $\widetilde{E}$ be a space
spanned by a basic sequence $\{ x_n\}_n$ such that $x_n\in E_n$
for $n\in\N$. Then no block sequence (with respect to $\{
x_n\}_n$) in $\widetilde{E}$ is $n-$sub\-symme\-tric for any
$n\in\N$, hence $\widetilde{E}$ contains no sub\-symme\-tric block
sequence, and by Lemma 3.3. no sub\-symme\-tric sequence.

Therefore by Theorem \ref{sub} the space $E_0$ contains a minimal
subspace. \hfill $\square$

\vspace{0.5cm}

Notice that we proved above that a Banach space saturated with
subsymmetric sequences  contains a "uniformly" minimal subspace,
ie. $c-$minimal for some $c\geq 1$.

\vspace{0.5cm}

\noindent {\it Proof of Theorem  \ref{sub}.}

We can assume that $E$ is a Banach space with a basis. We will use
below the notation introduced in the first section. Assume that
$E$ is saturated with $C-$sub\-symme\-tric sequences, for some
$C\geq 1$, and fix $\varepsilon >0$.

Pick $\delta >0$ satisfying $(C+\delta)^2(1+\delta)\leq
C^2+\varepsilon$. By Lemma \ref{bbasessub} and the density of $\q$
in $E$ the space $E$ is saturated with
$(C+\delta)-$sub\-symme\-tric block sequence from the family $\sgb
(E)$. We will prove that there is a block subspace $E_0\in\cgb
(E)$ such that every block subspace from the family $\cgb (E_0)$
contains a further block subspace $(C+\delta )^2-$iso\-morphic to
$E_0$. Therefore every infinitely dimensional subspace of $E_0$
contains a subspace $(C+\delta )^2(1+\delta )-$isomorphic to
$E_0$. By the choice of $\delta$ this will finish the proof of
Theorem
\ref{sub}.

From now on, unless otherwise stated, we consider block subspaces
from the family $\cgb (E)$ and sequences from the families
$\sgb (E)$ and $\sfb (E)$ only. Put $c=C+\delta$.

Recall that a tree $\ct$ on an arbitrary set $A$ is a subset of
the set $\bigcup_{n=1}^{\infty}A^n$ such that $\{ a_1,\dots
,a_n\}\in\ct$ whenever $\{ a_1,\dots ,a_n,a_{n+1}\}\in\ct$.

A branch of a tree $\ct$ is an infinite sequence $\{
a_n\}_{n\in\N}$ such that  $\{ a_1,\dots ,a_n\}\in\ct$ for any
$n\in\N$.

We will introduce now some notions. We call a tree
$\ct$ on $\q$ a block tree if $\ct\subset\sfb (E)$ and for any $\x\in\ct$ the set
$\ct (\x)=\{ x\in E:\; \x\cup \{ x \}\in\ct\}$ contains an infinite
block sequence in \q . Any branch of a block tree is a block
sequence. Moreover, since for any $\x
\in\ct$ we have $\ct (\x )\neq\emptyset$, every element $\x\in\ct$
is a part of some branch of $\ct$.

\bd Given sequences $\x, \y\in\sfb (E)$, $|\x|\geq |\y|$, a space
$L\in\cgb (E)$ and a tree $\ct$ on \q\ we write $(\x ;L)\sim(\y
;\ct )$ if $\ct=\bigcup\{ \ct_{\X}:\;\X\in\sgb (L),\; \X >\x\}$,
where $\{ \ct_{\X}\}$ are block trees on \q\ satisfying the
following conditions: \bnum\item for every block sequence
$\X\in\sgb (L)$, $\X >\x$ and every branch $\Y$ of $\ct_{\X}$ we
have $\Y >\y$ and sequences $\x\cup\X$, $\y\cup\Y$ are
$c-$equivalent,
\item for any block sequences $\X_1,\X_2\in\sgb (L)$, $\X_1 >\x$, $\X_2 >\x$ and
$n\in\N\cup\{ 0\}$, if $\X_1\cap E^n=\X_2\cap E^n$ then
$\ct_{\X_1}\cap E^{n+|\x|-|\y|}=\ct_{\X_2}\cap E^{n+|\x|-|\y|}$,
where $E^0=\{\Theta\}$.\enum\ed

This means that a tree of block sequences of $L$ beginning with a
finite sequence x can be represented in $\ct$ in a special manner.
In fact we will use the relation defined above only in the case
when $|\x |=|\y |$ or $|\x |=|\y |+1$.

\newpage

\bc\label{*} Take sequences $\x,\y\in\sfb (E)\cup\{\Theta\}$, $|\x
|\geq |\y |$,  a space $L\in\cgb (E)$, and a block tree $\ct$ on
\q . Assume $(\x;L)\sim (\y;\ct)$.
\bnum\item Let $x_0\in L$,
$\x<x_0$. Then there exists a block subtree $\ct'\subset\ct$ which
satisfies  $(\x\cup \{ x_0\};L)\sim (\y;\ct')$,
\item  Let $|\x |>|\y |$ and $y_0\in\ct \cap E$. Then
$\ct [y_0]=\{ \{ y_1,\dots,y_n\}:\; \{ y_0,y_1,\dots
,y_n\}\in\ct\}$ is a block tree and $(\x;L)\sim (\y\cup\{
y_0\};\ct [y_0])$.\enum\ec

\noindent {\it Proof of Claim \ref{*}.} For the first case, in
the situation as above define $\ct'$ by putting
$$\ct'_{\X}=\ct_{\{ x_0\}\cup\X}\;\;\; {\rm for}\;\;\;\X\in\sgb (L),\; \X >x_0.$$
The second case is obvious by the definition of the relation
$\sim$, since we put $(\ct [y_0])_{\X}=(\ct_{\X})[y_0]$ for any
$\X\in\sgb (L)$.\hfill $\square$

\vspace{0.5cm}

Now, given $M\in\cgb (E)$ put
$$\tau (M)=\left\{\begin{array}{c} (\x ,\y )\in (\sfb (E)\cup\{\Theta\})^2:
\; |\x |\geq |\y |,\; \exists\; L\in\cgb (E),\; L\leq M,\\
\exists\;\ct\; {\rm block\;\; tree\;\; \rm on}\;\; M\cap\q , \;
(\x;L)\sim (\y;\ct)\end{array}\right\}$$

Take $M_1\leq M_2$ and a pair $(\x,\y)\in\tau (M_1)$. Then there
is a space $L\leq M_1$ and a block tree $\ct_1$ on $M_1$ such that
$(\x;L)\sim (\y;\ct_1)$. Put $\ct_2=\ct\cap\bigcup_{n\in\N}
(M_2)^n$ (it means cutting off from $\ct_1$ sequences containing
vectors lying outside $M_2$). Then $\ct_2$ is also a block tree
(since $M_1\leq M_2$) satisfying  $(\x;L)\sim (\y; \ct_2 )$. One
only has to realize that for any sequence $\X\subset L$ a tree
$(\ct_2)_{\X}=\ct_{\X}\cap\bigcup_{n\in\N} (M_2)^n$ will do.

Therefore we have shown that the mapping $\tau$ is monotone, i.e.
if $M_1\leq M_2$ then $\tau (M_1)\subset\tau(M_2)$. Hence, on the
basis of Lemma \ref{stab}, there is a subspace $M_0\in\cgb (E)$
which is stabilizing for $\tau$.

\bc\label{**} Let $(\x,\y)\in\tau (M_0)$, $|\x|>|\y|$. Then for
any $M\in\cgb ( M_0)$ there is a vector $y_0\in M$ such that $(\x
,\y\cup\{ y_0\})\in \tau (M)$.\ec

\noindent {\it Proof of Claim \ref{**}.} In the situation as
above, by the stabilization property, for some subspace $L\in\cgb
(M)$ and a block tree $\ct$ on $M\cap\q$ we have $(\x ;L)\sim (\y
;\ct)$ and Claim \ref{*} finishes the proof of Claim \ref{**}.
\hfill $\square$

\vspace{0.5cm}

Given $M\in \cgb (M_0)$ put
$$\rho (M)=\left\{\begin{array}{c} (\x ,\y)\in (\sfb (E)\cup\{\Theta\})^2:\; |\x |\geq |\y |,\;
\exists \; L\in\cgb (E),\; L\doteq M, \\
\exists\;\ct \;{\rm block\;\; tree\;\; on}\;\; M_0\cap\q ,\;\;
(\x;L)\sim (\y;\ct)\end{array}\right\}$$ Take subspaces $M_1\geq
M_2$ and a pair $(\x,\y)\in\rho (M_1)$. There is a space $L_1$,
$L_1\doteq M_1$ and a tree $\ct_1$ on $M_0$ such that
$(\x;L_1)\sim (\y ; \ct_1)$. Put
$$L_2=L_1\cap M_2\doteq M_2,\;\;\;\;\ct_2=\bigcup \{(\ct_1)_{\X}:\;\X\in\sgb (L_2)\}$$
Then obviously
$(\x; L_2)\sim (\y ; \ct_2)$, hence $(\x,\y)\in\rho (M_2)$.

Therefore the mapping $\rho$ is monotone. Let $M_{00}\in\cgb
(M_0)$ be a stabilizing subspace for $\rho$, chosen on the basis
of Lemma \ref{stab}.

\bc\label{***} For any subspaces $M,N\in\cgb (M_{00})$ we have
$\rho (M)=\tau (N)$. \ec

\noindent {\it Proof of Claim \ref{***}.} By the stabilization
property it is enough to prove that $\tau (M_{00})=\rho (M_{00})$.
By definition and the stabilization property $\rho
(M_{00})\subset\tau (M_0)=\tau (M_{00})$. Now, if $(\x, \y)\in\tau
(M_{00})$, then $(\x,\y)\in\rho (M)$ for some $M\subset M_{00}$,
hence again by the stabilization $(\x,\y)\in\rho (M_{00})$. \hfill
$\square$

\vspace{0.5cm}

By the assumption and Lemma \ref{bbasessub} there is a
$c-$sub\-symme\-tric sequence $\{ z_n\}_{n=1}^{\infty}\in\sgb
(M_{00})$. Let $E_0=\spa\{z_n\}_n$.

\bc\label{****} $(\{\Theta\} ;E_0)\sim (\{\Theta\} ; \sfb (E_0))$,
in particular $(\{ \Theta\} ,\{\Theta\} )\in\tau (M_{00})$. \ec

\noindent {\it Proof of Claim \ref{****}.} Take any block
sequence $\X=\{ x_n\}_{n=1}^{\infty}\in\sgb (E_0)$. Then
$$x_n=\sum_{i=i_{n}}^{i_{n+1}-1}a_i\; z_i, \;\;\;\; n\in\N$$
for some scalars $\{ a_i\}_{i\in\N}$ and some sequence $\{i_n\}_{n\in\N}\subset\N$. Put
$$\ct_{\X}=\left\{\;\;\left\{\sum_{i=i_{1}}^{i_{2}-1}a_i\; z_{\phi (i)},\dots ,\sum_{i=i_{n}}^{i_{n+1}-1}a_i\; z_{\phi(i)}
\right\}: \begin{array}{c}\phi :\N\rightarrow\N\; {\rm strictly\;\; increasing},\\ n\in\N\\
\end{array}\right\}$$
Obviously every set $\ct_{\X}$ is a block tree. Moreover, $\sfb
(E_0)=\bigcup\{\ct_{\X}:\;\X\in\sgb (E_0)\}$ and, by
$c-$subsymmetry of the sequence $\{ z_n\}_n$, for any $\X\in\sgb
(E_0)$ every infinite branch of $\ct_{\X}$ is $c-$equivalent to X.
The "uniqueness" condition is also satisfied. \hfill $\square$

\vspace{0.5cm}

We will show that every block subspace from the family  $\cgb
(E_0)$ contains a further subspace $c^2-$iso\-morphic to $E_0$
which will finish the proof of the Theorem.

Take arbitrary $M\in\cgb (E_0)$. We will pick by induction block
sequences $\{ z_{k_n}\}$ and $\{ y_n\}\subset M$ such that $(\z_n
,\y_n)\in\tau (M)$ for $n\in\N$, where $\z_n=\{ z_{k_1},\dots,
z_{k_n}\}$ and $\y_n=\{ y_1,\dots ,y_n\}$, $n\in\N$. By definition
it implies in particular that for any $n\in\N$ sequences $\{
z_{k_1},\dots, z_{k_n}\}$ and $\{ y_1, \dots ,y_n\}$ are
$c-$equivalent, thus also sequences $\{ z_{k_n}\}_{n\in\N}$ and
$\{ y_n\}_{n\in\N}$ are $c-$equivalent. By $c-$subsymmetry of the
sequence $\{ z_n\}_{n\in\N}$ sequences $\{ z_n\}_n$ and $\{
y_n\}_n$ are $c^2-$equivalent, hence $E_0$ shares the demanded
property.

Put $k_1=1$. By Claims \ref{****} and \ref{*}
$(\z_1,\{\Theta\})\in\tau (M_0)=\tau (M)$. By Claim \ref{**}
there is a vector $y_1\in M$ such that $(\z_1,\y_1)\in\tau (M)$.

Assume now that we have picked vectors $z_{k_1},\dots,z_{k_n}$ and
$y_1,\dots, y_n\in M$ such that $(\z_n,\y_n)\in\tau (M)$. By
Claim \ref{***} $(\z_n,\y_n)\in\rho (E_0)$. Therefore for some
$L\doteq E_0$ there is a tree $\ct$ on $M_0$ such that $(\z_n
;L)\sim (\y_n;\ct)$. Let $k_{n+1}$ be such that
$z_{k_{n+1}}>z_{k_n}$ and $z_{k_{n+1}}\in L$. Then by Claim
\ref{*} $(\z_{n+1},\y_n)\in\tau (M_0)$. Hence by Claim \ref{**}
there is a vector $y_{n+1}\in M$ such that
$(\z_{n+1},\y_{n+1})\in\tau (M)$, which finishes the inductive
step and the proof of Theorem \ref{sub}. \hfill $\square$

\vspace{0.5cm}

Obviously in the reasoning above we did not use the
unconditionality property, never\-theless due to Gowers' dichotomy
concerning HI spaces and unconditional sequences as well as
properties of HI spaces (the lack of non-trivial isomorphisms) we
can assume without the loss of generality that we are dealing with
spaces with unconditional bases.

\end{document}